\documentclass{amsart}
\usepackage{amsmath,amsfonts}
\usepackage{amscd,amsthm}

\theoremstyle{plain} 

\newtheorem{theorem}{Theorem}[section]
\newtheorem{corollary}[theorem]{Corollary}
\newtheorem{lemma}[theorem]{Lemma}
\newtheorem{proposition}[theorem]{Proposition}

\newtheorem{remark}[theorem]{Remark}

\newcommand{\R}{\mathop{\mathbb{R}}}
\newcommand{\Q}{\mathop{\mathbb{Q}}}
\newcommand{\Z}{\mathop{\mathbb{Z}}}
\newcommand{\N}{\mathop{\mathbb{N}}}
\newcommand{\C}{\mathop{\mathbb{C}}}

\newcommand{\diag}[1]{\mathop{\it diag(#1_1,\dots,#1_n)}}

\numberwithin{equation}{section}

\begin{document}

\title[Nekhoroshev--like estimate for non--linearizable 
analytic vector fields.]{Exponentially long time stability 
near an equilibrium point
for non--linearizable analytic vector fields.}

\author{Timoteo Carletti}

\date{\today}

\address[Timoteo Carletti]{Scuola Normale Superiore piazza dei Cavalieri 7,
 56126  Pisa, Italia}

\email[Timoteo Carletti]{t.carletti@sns.it}



\begin{abstract}
We study the orbit behavior of a germ of an analytic vector field of
$(\C^n,0)$, $n \geq 2$. We prove that if its linear part is semisimple,
non--resonant and verifies a Bruno--like condition, then the origin is
effectively stable: stable for finite but exponentially long times.
\end{abstract}

\maketitle

\section{Introduction}
Let us consider the germ of analytic vector field, $X_F=\sum_{1\leq j 
  \leq n}F_j(z)\frac{\partial}{\partial z_j}$, of $(\C^n,0)$ $n \geq 2$, whose
  components $(F_j)_{1\leq j \leq n}$ are analytic functions vanishing 
  at $0\in \C^n$. 

Let us consider the associated Ordinary Differential Equation:
\begin{equation}
\label{eq:ODE}
\frac{dz}{dt}=F(z) \, ;
\end{equation}
under the above assumptions $z(t;0)= 0$ for all $t$ is an 
equilibrium solution~\footnote{Here and throughout the paper by
$z(t;z_0)$ we mean the solution at time $t$ of~\eqref{eq:ODE}
  s.t. $z(0;z_0)=z_0$. When the value of $z_0$ will not be relevant
  we'll just write $z(t)$.}. We are interested in studying the stability of  
orbits of $X_F$ in a neighborhood of this equilibrium point.

We use the standard definition of {\em stability} (see~\cite{Moser}) for
an equilibrium solution: $z=0$ it is stable is the past and in the future if
for any neighborhood $U$ of $0$ there exists a neighborhood V, containing the
origin, s.t. $z(0;z_0)\in V$ implies $z(t;z_0)\in U$ for all $t\in \R$.

In a coordinates system
 centered at the equilibrium point the $j$--th component of the vector field 
will take the form: $F_j(z)=(Az)_j+f_j(z)$, with $A$ a $n\times n$
 complex matrix and 
 $f_j$ analytic function such that $f_j(0)=Df_j(0)=0$, for all $1\leq
 j \leq n$. 
Following the idea of Poincar\'e to study the orbit of~\eqref{eq:ODE} in a 
neighborhood of the origin, we will try to find an analytic change of
 coordinates, 
through an analytic diffeomorphisms $z\mapsto H(z)=w$ the {\em linearization},
 s.t. in the new coordinates the vector field $X_F$ is conjugate to
 its linear part, $X_A=\sum (Az)_j \frac{\partial}{\partial z_j}$:
 $H^* X_F H^{-1}=X_A$. Hence equation~\eqref{eq:ODE} rewrites:
\begin{equation}
\label{eq:ODE2}
\frac{dw}{dt}=Aw \, .
\end{equation}
This change of coordinates must solve:
\begin{equation}
AH(z)=DH(z) \cdot \left( Az+f(z) \right) \, ,
\label{eq:conjugacy}
\end{equation}
and it is unique by assuming $DH(0)=\mathbb{I}$.

Clearly if the linear system~\eqref{eq:ODE2} is stable and~\eqref{eq:ODE} is 
analytically linearizable, then also the latter is stable. It is a remarkable result
that this condition is also necessary, as the following Theorem states:

\begin{theorem}[Carath\'eodory--Cartan 1932]
\label{the:caratheodorycartan}
Necessary and sufficient condition for the stability of the solution $z=0$ 
of~\eqref{eq:ODE} for all real $t$ is that:
\begin{enumerate}
\item $A$ is diagonalizable with purely imaginary eigenvalues;
\item there exists an holomorphic function
  $z=K(w)=w+\mathcal{O}(|w|^2)$, $w\in\C^n$, which 
brings~\eqref{eq:ODE} into the linear system:
\begin{equation*}
\frac{dw}{dt}=Aw \, .
\end{equation*}
\end{enumerate}
\end{theorem}

So let us assume $A$ to verify hypothesis of
Theorem~\ref{the:caratheodorycartan}: 
let $(\omega_j)_{1\leq j \leq n}\subset \R$ and $A=\diag{i\omega}$.

Then $A$ belongs to the {\em Siegel domain}~\footnote{According to the 
classification of~\cite{Bruno} this case is {\em Poincar\'e domain 1.d},
 but we prefer consider it as a Siegel case because the obstructions to the
linearizability are very similar to those encountered in the Siegel domain.}:
 the origin is contained in the convex hull of the set  
of eigenvalues plotted as points in the complex plane (e segment in
this case). This is the harder situation
w.r.t. to the complementary case, {\em Poincar\'e domain},
 because {\em small divisors}
are involved: the existence of an analytic linearization is strictly related 
to the arithmetic property of approximation of the vector
$\omega=(\omega_1,\dots,\omega_n)$, with vectors with integer entries.

The first step is to assume $A$ to be {\em non--resonant}: 
$\alpha \cdot \omega \neq \omega_j$, for all $\alpha \in \N^n$ s.t. 
$|\alpha|=\alpha_1+\dots +\alpha_n \geq 2$ and for 
all $j\in \{ 1, \dots, n\}$. This ensures the existence
of a {\em formal} change of variable which linearizes~\eqref{eq:ODE}.

In~\cite{Bruno} author introduced the, today called, {\em Bruno
  condition}~\footnote{ 
The Bruno condition can be rewritten using a general increasing
  sequence of integer  
numbers, $(p_k)_k$. In~\cite{Bruno} pag. 222, author proved
  that~\eqref{eq:brunovf} is equivalent to:
\begin{equation*}
\sum_{k \geq 0}\frac{\log \Hat{\Omega}^{-1}(p_{k+1})}{p_{k}}<+\infty \, ,
\end{equation*}
where $\Hat{\Omega}(p)=
\min\{ |\alpha \cdot \omega - \omega_j|: j\in \{ 1,\dots, n\}, \alpha \in\Z^n,
 0<|\alpha |<p,\}$.} to
characterize the rate 
of approximation of vectors of $\R^n$ with vectors of $\Z^n$:
\begin{equation}
\label{eq:brunovf}
\sum_{k \geq 0}\frac{\log \Omega^{-1}({k+1})}{2^{k}}<+\infty \, ,
\end{equation}
where for all positive integers $k$, $\Omega(k)=
\min\{ |\alpha \cdot \omega - \omega_j|: j\in \{ 1,\dots, n\}, \alpha \in\Z^n,
 0<|\alpha |<2^k\}$. 
From the result of Bruno~\cite{Bruno} follows that if $A$ is
non--resonant, diagonal with purely  
imaginary eigenvalues and verifies the above condition, then there exists an 
analytic linearization which brings~\eqref{eq:ODE} into~\eqref{eq:ODE2}.

Let us make a step backward and consider the following problem~\cite{Carletti2003}.
Let $\mathcal{A}_1 \subset \mathcal{A}_2 \subset 
\mathbb{C}^n \left[ \left[ z_1,\dots ,z_n \right] \right]$ be two classes of
 formal power series closed w.r.t. to derivation and composition. Let 
$\hat{f}\in \mathcal{A}_1$ s.t. $\hat{f}=\sum_{|\alpha|\geq
   2}f_{\alpha}z^{\alpha}$, let $A\in GL(n,\C)$ and consider 
 the following (formal) ODE:
\begin{equation}
\label{eq:ODEA1}
\frac{d z}{dt}=Az+\hat{f}(z) \, .
\end{equation}
We say that~\eqref{eq:ODEA1} is {\em linearizable in $\mathcal{A}_2$} 
if there exists $\hat{H}\in \mathcal{A}_2$, normalized with
 $\hat{H}=z+\mathcal{O}(|z|^{|\alpha|})$, 
$|\alpha| \geq 2$, s.t. {\em formally} we have:
\begin{equation*}
w=\hat{H}(z) \quad \text{and} \quad \frac{d w}{dt}=Aw \, .
\end{equation*}

If both $\mathcal{A}_1$ and $\mathcal{A}_2$ coincide with the ring of formal
 power series we already know that generically the problem has
 solution if and only if $A$ 
 is {\em non--resonant}, which will be assumed from now. In the other cases of
general algebras of formal power series, new arithmetical conditions on $A$
 have to be imposed if we are in the Siegel domain. This case has been considered
in details in~\cite{Carletti2003} section $5$, to which we refer for all 
details~\footnote{See also~\cite{CarlettiMarmi2000} where a similar
problem for germs of diffeomorphisms of $(\C,0)$ has been
 studied.}. There author 
proved that the Bruno condition is still sufficient to linearize whenever
$\mathcal{A}_1=\mathcal{A}_2$, otherwise new Bruno--like conditions are introduced,
 weaker than the original Bruno condition.

An interesting case is when $\mathcal{A}_1$ is the ring of convergent power
series in some neighborhood of the origin, and $\mathcal{A}_2$ is the algebra of
 {\em Gevrey}--$s$, $s>0$, formal power series. Namely we are
considering the {\em Gevrey linearization of analytic vector fields}. 

 Let $\Hat F=\sum f_{\alpha} z^{\alpha}$, 
$(f_{\alpha})_{\alpha \in \N^n} \subset \C^n$ be a 
formal power series, then we say that it is 
{\em Gevrey--$s$}~\cite{Balser1994,Ramis1991},
 $s>0$, if there exist two positive constants $C_1,C_2$ such that:
\begin{equation}
  \label{eq:gevreydefvect}
  |f_{\alpha}| \leq C_1 C_2^{-s|\alpha|} (|\alpha|!)^s  \quad \forall
   \alpha \in \mathbb{N}^n \, . 
\end{equation}
We denote the class of formal vector valued power series Gevrey--$s$ 
by $\mathcal{C}_s$. It is closed w.r.t. derivation and composition.

In the Gevrey--$s$ case the arithmetical condition introduced 
in~\cite{Carletti2003}, called {\em Bruno}--$s$ condition, $s>0$, 
 for short $\mathcal{B}_s$, reads:
\begin{equation}
  \label{eq:brunosndim}
  \limsup_{|\alpha|\rightarrow +\infty}\left( 2\sum_{m=0}^{\kappa(\alpha)} 
\frac{\log \Omega^{-1}(p_{m+1})}{p_m}-s\log |\alpha|\right)< +\infty \, ,
\end{equation}
for some increasing  sequence of positive integer $(p_k)_k$ and
$\kappa(\alpha)$ is defined by $p_{\kappa(\alpha)}\leq |\alpha| < p_{\kappa(\alpha)+1}$. 

\begin{remark}
  This definition recall the classical one of
   Bruno~\cite{Bruno}, where first one suppose the existence of a
   strictly increasing sequence of positive integer such
   that~\eqref{eq:brunosndim} holds, then one can prove (see~\cite{Bruno}
   \S IV page 222) that one can take an exponentially growing
   sequence, e.g. $p_k=2^k$. This holds also in our case, in fact we
   can prove that~\eqref{eq:brunosndim} is equivalent to:  
\begin{equation*}
  \limsup_{N\rightarrow +\infty}\left(\sum_{l=0}^N\frac{\log
  \Omega^{-1}(2^{l+1})}{2^l}-sN2\log2 \right) < +\infty \, .
\end{equation*}
A proof of this claim can be found in~\cite{Carletti2002}.
\end{remark}

When $n=2$, under the above condition (non--resonance and Siegel domain),
rescaling time by $-\omega_2$ (assuming $\omega_2 \neq 0$),
 the ODE associated to the vector field can be rewritten as:
\begin{equation}
\label{eq:ODEn2}
\begin{cases}
\dot z_1 = \omega z_1 + h.o.t. \\
\dot z_2 = - z_2 + h.o.t. \\
\end{cases}\, ,
\end{equation}
where $\omega = -\omega_1/\omega_2\in (\R\setminus\Q)^{+}$ and high
order terms  
means $\mathcal{O}(|z|^{|\alpha|})$ with $|\alpha|\geq 2$, namely only
the ratio of the eigenvalues enters. Then the Bruno--$s$ condition
 can be slightly weakened (see~\cite{CarlettiMarmi2000}):
\begin{equation}
\label{eq:brunos1dim}
\limsup_{n\rightarrow +\infty} \left( \sum_{j=0}^{k(n)}\frac{\log
  q_{j+1}}{q_{j}} - s\log n \right) <+\infty \, , 
\end{equation}
where $k(n)$ is defined by $q_{k(n)}\leq n < q_{k(n)+1}$ and $(q_n)_n$ are
the denominators of the convergents~\cite{HardyWright} to $\omega$. 

We remark that
 in both cases the new conditions are weaker than Bruno's condition, which
is recovered when $s=0$. When $n=2$ we prove that the set $\bigcup_s \mathcal{B}_s$ 
is $PSL(2,\Z)$--invariant (see remark~\ref{rem:invariance}). 

The main result of~\cite{Carletti2003} in the case of Gevrey--$s$ classes reads:
\begin{theorem}[Gevrey--$s$ linearization]
\label{thm:gevreylin}
Let $\omega_1,\dots,\omega_n$ be real numbers and $A=\diag{i\omega}$;
let $D_1 = \{ z \in \C^n : |z_i|<1 \, , 1\leq i\leq n \}$ be the isotropic polydisk 
of radius $1$ and let $F:D_1\rightarrow \C^n$ be an analytic function,
such that $F(z)=Az+f(z)$, with $f(0)=Df(0)=0$. If $A$ is non--resonant and verifies
a Bruno--$s$, $s>0$, condition~\eqref{eq:brunosndim} (or
condition~\eqref{eq:brunos1dim} if $n=2$),
then there exists a formal Gevrey--$s$ linearization $\Hat{H}$.
\end{theorem}
The aim of this paper is to show that the Gevrey character of the formal 
linearization can give information concerning the dynamics of the analytic vector
field. Let $F(z)=Az+f(z)$ verify hypotheses
 of Theorem~\ref{thm:gevreylin}, assume moreover $X_F$ not to be analytically 
linearizable. We will show that even if there is not a {\em Stable domain}, 
where the dynamics of $X_F$ is conjugate to the dynamics of its linear
part, we have  
an open neighborhood of the origin which ``behaves as a Stable domain'' for the
flow of $X_F$ for finite but long time, which results exponentially long:
the {\em effective stability}~\cite{GFGS,GiorgilliPosilicano} 
of the equilibrium solution.

In the case of analytic linearization, $|H_j(z)|$, $j=1, \dots, n$,
is {\em constant along the orbits}, namely it is a {\em first
  integral} and the flow of~\eqref{eq:ODE} is bounded 
 for all $t$ and sufficiently small $|z_0|$.

We will prove that any non--zero $z_0$ belonging to a polydisk of sufficiently
 small radius $r>0$, can be followed up to a time 
$T=\mathcal{O}(exp \{ r^{-1/s} \} )$, being $s>0$ the Gevrey
exponent of the formal linearization, and we can find an 
{\em almost first integral}: 
a function which varies by a quantity of order $r$ during this
 interval of time.
More precisely we prove the following

\begin{theorem}
  \label{thm:maintheorem}
Let $n\in\N$, $n\geq 2$. Given real $\omega_1,\dots,\omega_n$ consider 
$A=\diag{i\omega}$; let $F:D_1\rightarrow \C^n$ be an analytic function,
such that $F(z)=Az+f(z)$, with $f(0)=Df(0)=0$. If $A$ is non--resonant
and verifies 
a Bruno--$s$, $s>0$, condition~\eqref{eq:brunosndim}
(or~\eqref{eq:brunos1dim} if 
$n=2$),
then for all sufficiently small $0< r_{**} <1$, there exist positive constants 
$A_{**},B_{**},C_{**}$ such that for all $0<|z_0|<r_{**}/2$, the
solutions $z(t;z_0)$ 
are well defined and  verify $|z(t;z_0)|\leq C_{**} r_{**}$, for all 
$|t|\leq T_* = A_{**}^{-1}\, exp 
\Big \{ B_{**}\left( r_{**}/|z_0| \right)^{1/s} \Big \}$.
\end{theorem}

We want to stress here that when $s\rightarrow 0$ the stability time
  goes to infinity, because the {\em critical exponent of stability
  time} is $1/s$: namely one get {\em stability}. At the same time the
  Bruno--$s$ condition ``tends'' to the classical Bruno condition,
  which is a sufficient condition to ensure analytic linearizability
  and hence stability under our assumptions. So we recover the classical stability
result as limit of longer and longer effective stability times.

Results similar to Theorem~\ref{thm:maintheorem} have been obtained
in~\cite{GFGS} for hamiltonian 
vector fields and in~\cite{GiorgilliPosilicano} for real reversible systems
of coupled harmonic oscillators. In both papers effective stability is
proved by assuming the linear part of the vector field to verify some
Diophantine condition~\footnote{A vector $\omega \in \R^n$ belongs to
  $CD(\gamma,\tau)$ if there exists $\gamma>0$ and $\tau > n-1$ such
  that for all $\alpha \in \N^n$ and all $j\in \{1,\dots,n\}$ one has:
  $|\alpha \cdot \omega - \omega_j| \geq \gamma |\alpha|^{-\tau}$. Let
  $A=\diag{\omega}$, then $A$ verifies a Diophantine condition if
  $\omega$ does.} $CD(\gamma,\tau)$, for some $\gamma>0$ and $\tau >
n-1$, and the critical exponent of stability time is $1/\tau$. In our
result, too, the critical exponent of stability time depends on some
arithmetical property of the linear part of the vector field but in a
more general way in fact we assume $A \in \mathcal{B}_s \supset
CD(\gamma,\tau)$, for all $\gamma >0$ and $\tau \geq n-1$.

The second remark is that in~\cite{GFGS,GiorgilliPosilicano} effective
stability is obtained using some partial normal form, then working
on it and using the Poincar\'e summation at the smallest term 
(see Lemma~\ref{lem:sumupsmallest}), the proof is done.
Here the method used is completely different: we first
linearize formally the system and then using properties of the formal
linearization we conclude still using the Poincar\'e summation at the smallest term.
 This method introduce also our main
drawback: we must assume $A$ to be non--resonant (to linearize) and
this prevents us from considering real vector fields and hamiltonian
ones, where an ``intrinsic'' resonance is present.

In section~\ref{sec:conclusions} we discuss the relation between the
Bruno--$s$ condition  and other arithmetical conditions.


\section{Proof of the main Theorem}
\label{sec:proofmainthm}

In this part we will prove our main result, Theorem~\ref{thm:maintheorem}. The
proof will be divided into three steps: first we use the Gevrey--$s$ 
character of the formal linearization $\Hat H$, given by 
Theorem~\ref{thm:gevreylin}, to find an approximate solution of the 
conjugacy equation~\eqref{eq:conjugacy} up
to a (exponentially) small correction (paragraph~\ref{ssec:firststep});
then we prove a Lemma allowing us to control how the small 
error introduced in the solution propagates (paragraph~\ref{ssec:thirdstep}). 
Finally we collect all
the informations to conclude the proof (paragraph~\ref{ssec:endproof}).

\subsection{Determination of an approximate solution}
\label{ssec:firststep}

Let $F$ verifies hypotheses of Theorem~\ref{thm:maintheorem} and let us consider the
 first order differential equation in $\mathbb{C}^n$, 
 $n \geq 2$:
\begin{equation}
\label{eq:ode}
\frac{dz}{dt} = F(z) \, . 
\end{equation}
By Theorem~\ref{thm:gevreylin} this system can be put in linear form
by a formal power series  
$\Hat{H}$ which belongs to $\mathcal{C}_s$ and it solves (formally):
\begin{equation}
  \label{eq:forh}
  \frac{d}{dt}\Hat H(z) = A \Hat H(z) \, ,
\end{equation}
we observe that one can choose $\Hat H(z)=z+\mathcal{O}(|z|^2)$.

Since ${\Hat H}=\sum h_{\alpha} z^{\alpha}\in \mathcal{C}_s$,
 there exist positive 
constants $A_1$ and $B_1$ such that 
\begin{equation}
\label{eq:gevreyhm1}
|h_{\alpha}| \leq A_1 B_1^{-s|\alpha|} (|\alpha|!)^s 
\quad \forall \,|\alpha| \geq 1 \, .
\end{equation}
For any positive integer $N$ we consider the {\em vectorial polynomial}, 
sum of homogeneous vector monomials of degree $1\leq l \leq N$,
defined by: $\mathcal{H}_N(z)=\sum_{l=1}^{N} \sum_{|\alpha|=l}
 h_{\alpha} z^{\alpha}$ and the {\em Remainder Function}:
\begin{equation}
\label{eq:remainderfunction}
\mathcal{R}_N(z)=D \mathcal{H}_N (z)\cdot F(z) 
- A\mathcal{H}_N(z) \, .
\end{equation}
Clearly $\mathcal{H}_N(z)$ doesn't solve the linearization problem, but:
\begin{equation}
\label{eq:forhN}
\frac{d}{dt}\mathcal{H}_N(z) = A \mathcal{H}_N(z)+ \mathcal{R}_N(z)\, ,
\end{equation}
hence the remainder function gives the difference from the true solution 
and the approximate one.

The following Proposition collects some properties of the remainder 
function.

\begin{proposition}
Let $\mathcal{R}_N(z)$ be the remainder function defined 
in~\eqref{eq:remainderfunction} and
 let $\alpha \in \N^n$, then:
  \begin{enumerate}
  \item[1)] $\partial_z^{\alpha} \mathcal{R}_N(0)=0$ if $|\alpha| \leq N$.
  \item[2)] For all $0<r<1$ there exist positive constants $A_2$ and 
$B_2$
 such that if $|\alpha| \geq  N+1$, then:
\begin{equation*}
\Big |\frac{1}{\alpha !}\partial_z^{\alpha} 
\mathcal{R}_N(0)\Big |\leq A_2 
r^{-|\alpha|}B_2^{-sN} (N!)^s\, .
\end{equation*}
  \item[3)] For all $0<r<1$ and $|z|<r/4$ there exist positive 
constants $A_3,B_3$ such that:
\begin{equation}
\label{eq:punto3}
 |\mathcal{R}_N(z)| \leq A_3 B_3^{-sN} (N!)^s
\left(\frac{|z|}{r}\right)^{N+1} \, .
\end{equation}
\end{enumerate}
Where we used the compact notation 
$\frac{1}{\alpha !}\partial_z^{\alpha} =
\frac{1}{\alpha_1 ! \dots \alpha_n !}
\frac{\partial^{|\alpha|}}{\partial_{z_1}^{\alpha_1}
\dots \partial_{z_n}^{\alpha_n}}$.
\end{proposition}

\proof 
Statement 1) is an immediate consequence of 
the definition of $\mathcal{R}_N$.

To prove 2) we observe that $\mathcal{R}_N(z)$ is an analytic function 
on $D_1$, being obtained with product of analytic functions, then one 
gets by Cauchy's estimates for all $0<r<1$ and for all $|\alpha|\geq N+1$:
\begin{equation}
  \label{eq:cauchy}
\Big | \frac{1}{\alpha!}\partial_z^{\alpha}\mathcal{R}_N(0)\Big |
\leq \frac{1}{(2\pi)^n} 
\frac{1}{r^{|\alpha|+1}}\max_{|z|=r}|D \mathcal{H}_N\cdot F(z)| \, ,
\end{equation} 
because $\partial^{\alpha}\mathcal{H}_N=0$ for $|\alpha|\geq N+1$.
Recalling the Gevrey estimate~\eqref{eq:gevreyhm1} for $\mathcal{H}_N$ 
and the analyticity of $F$ we obtain:
\begin{equation}
  \label{eq:estim2}
\Big | \frac{1}{\alpha!}\partial_z^{\alpha}\mathcal{R}_N(0)\Big |
\leq A_2 B_2^{-sN} (N!)^s r^{-|\alpha|} \, ,  
\end{equation}
for some positive constants $A_2$ and $B_2$ depending on the previous 
constants, on the dimension $n$ and on $F$.

To prove 3) let us write the Taylor series $\mathcal{R}_N(z)=
\sum_{|\alpha|\geq N+1}\frac{1}{\alpha!}\partial_z^{\alpha}
\mathcal{R}_N(0) z^{\alpha}$: 
the bound on derivatives~\eqref{eq:estim2} implies
the estimate~\eqref{eq:punto3} for all $|z|< r/4$ and for some 
positive constants $A_3$ and $B_3$. 
\endproof

The bound~\eqref{eq:punto3} on $\mathcal{R}_N(z)$ depends on the positive 
integer $N$, so we can determine the value of $N$ for which the right
hand side of~\eqref{eq:punto3} attains its minimum, that's 
Poincar\'e's idea
of {\em summation at the smallest term}.

\begin{lemma}[Summation at the smallest term]
 \label{lem:sumupsmallest}

Let $\mathcal{R}_N(z)$ defined as before
and let $0<r_*<1/4$ then there exist positive constants $A_4,B_4$
 such that for all $0<|z|<r_*$ we have:
\begin{equation}
  \label{eq:boundsumup}
  |\mathcal{R}_{\bar{N}}(z)|\leq A_4 \, exp \Big \{ -B_4
\left(\frac{r_*}{|z|}\right)^{1/s} \Big \} \, ,
\end{equation}
where $\bar{N}=\lfloor B_4\left(r_*/|z|\right)^{1/s}\rfloor$ and 
$\lfloor x \rfloor$ denotes the integer part of $x\in \R$.
\end{lemma}

\proof
Let us fix $0<r_*<1/4$, then for $0<|z|<r_*$ by Stirling formula
we obtain:
\begin{equation}
\label{eq:lem1}
| \mathcal{R}_N (z) | \leq A_4 \left( N B_3^{-1} \left(|z|/r_*\right)^{1/s} \right)^{Ns} e^{-sN} \, ,
\end{equation}
for some positive constant $A_4$. The right hand side of~\eqref{eq:lem1} attains its 
minimum at $\bar{N} = B_3 \left(r_*/|z|\right)^{1/s}$, evaluating the value of this minimum 
we get~\eqref{eq:boundsumup} with $B_4=B_3$.
\endproof

\subsection{Control of the ``errors''}
\label{ssec:thirdstep}

Let us define $\mathcal{H}(z)=\mathcal{H}_{\bar{N}}(z)$ and $\mathcal{R}(z)=
\mathcal{R}_{\bar{N}}(z)$, being $\bar{N}$ the ``optimal value'' obtained in
Lemma~\ref{lem:sumupsmallest}. We remark that $\mathcal{H}(z)$ doesn't 
solve~\eqref{eq:forh} but the ``error'', $\mathcal{R}(z)$, can be made
very small: exponentially small. We will prove that for initial 
conditions in a sufficiently small disk, one can follow the flows 
for an exponential long time without leaving a disk comparable size.

\begin{lemma}[Control the flow]
\label{lem:iteration}
Let $a,b,\alpha$ and $R$ be positive real numbers. Let 
$T=Ra^{-1}e^{b/(2R)^{\alpha}}$ and let us consider the 
Cauchy problem:
\begin{equation*}
\begin{cases}
\frac{d}{dt}x(t)=a e^{-b/x^{\alpha}} \\
x(0)=R \quad .
\end{cases}
\end{equation*}
Then $0<x(t)<2R$ for all $|t|< T$.
\end{lemma}
\proof
Let us write the Cauchy problem in integral form:
\begin{equation*}
x(t)=R+\int_0^t a e^{-b/(x(s))^{\alpha}} \, ds \, ,
\end{equation*}
$x(t)$ is trivially monotonically increasing, hence the same holds
for the function $t\mapsto a e^{-b/(x(t))^{\alpha}}$.
Let us suppose that there exists $0<t_0 <T$, for which $x(t_0)=2R$; then
\begin{equation*}
2R=x(t_0)=R+\int_0^{t_0} a e^{-b/(x(s))^{\alpha}} \, ds <R
+t_0 a e^{-b/(x(t_0))^{\alpha}} \, ,
\end{equation*}
namely $t_0>Ra^{-1} e^{b/(2R)^{\alpha}}=T$, which gives a contradiction.
Hence either $x(t_0)>2R$ for all $0< t <T$ or $x(t_0)<2R$, but the first
case have to be excluded because $x(0)=R<2R$.

The case $t<0$ can be handle in a similar way by showing that $t\mapsto x(t)$ doesn't
decrease too much.
\endproof

Let $r_*$ as in Lemma~\ref{lem:sumupsmallest}, define 
$\rho(z)=|\mathcal{H}(z)|$ for 
all $0<|z|<r_*$, then Lemma~\ref{lem:sumupsmallest} admits the following 
Corollary, which allows us to control the evolution of $\rho(z)$.

\begin{corollary}
  \label{lem:itarate}
Let $0<r_*<1/4$, then there exists $0<r_{**}\leq r_*$ and 
positive constants $A_*,B_*$ such that for all $0<|z|<r_{**}$
 we have:
\begin{equation}
\label{eq:iterate}
\frac{d}{dt} \rho(z)\leq A_* \, exp \Big \{ -B_*
\left( \frac{r_*}{\rho(z)} \right)^{1/s} \Big \} \, .
\end{equation}
\end{corollary}

\proof
Let $j\in \{ 1, \dots ,n\}$ and let us consider the time evolution of
$|\mathcal{H}_j(z(t))|$. If $\mathcal{H}$ was a solution this would be
a constant of motion, this is not the case but its evolution is
nevertheless very slow. In fact thanks to~\eqref{eq:forhN} we get:
\begin{equation*}
\frac{d}{dt}|\mathcal{H}_j(z)|\leq |\mathcal{R}_j(z)| \, ,
\end{equation*}
hence a similar statement holds for $\rho(z)=|\mathcal{H}(z)|=
\max_{1\leq j \leq n}|\mathcal{H}_j(z)|$.

We want now to express the
exponential smallness of $|\mathcal{R}(z)|$ in terms of $\rho(z)$
 instead of $|z|$. $\mathcal{H}(z)$ is tangent to the identity close 
to zero and then locally invertible. The inverse is still tangent to
the identity, vanishing at zero and analytic in a neighborhood of the
origin, then sufficiently close
to the origin we have $|\mathcal{H}^{-1}(w)|\leq C|w|$, for some $C>0$.
Finally we can take $|z|$ sufficiently small, say $|z|<r_{**}$ for some 
$0<r_{**}\leq r_*$, s.t. $|z|\leq C |\mathcal{H}(z)|<r_*$, hence:
\begin{equation*}
|\mathcal{R}(z)|\leq A_4 \, exp \Big \{ -B_4 
\left(\frac{r_*}{|z|}\right)^{1/s} \Big \} \leq 
A_4 \, exp \Big \{ -B_4 C^{-1/s} 
\left(\frac{r_*}{\rho(z)}\right)^{1/s} \Big \} \, ,
\end{equation*}
and the claim follows with $A_*=A_4$ and $B_*=B_4C^{-1/s}$.
\endproof

\subsection{End of the proof}
\label{ssec:endproof}

We are now able to conclude the proof of the main 
Theorem~\ref{thm:maintheorem}. Take any $0<|z_0|<r_{**}/2$ and let 
$\rho_0=|\mathcal{H}(z(0;z_0))|$, then there exists a positive constant $C_1$
s.t. $\rho_0\leq C_1|z_0|$. By Corollary~\ref{lem:itarate} we have
\begin{equation}
  \label{eq:firstest}
\frac{d}{dt}  \rho(z(t;z_0)) \leq A_* \, exp \Big \{ -B_* 
\left( r_*/\rho(z(t;z_0))\right)^{1/s} \Big \} \, .
\end{equation}
Let us call $R=\rho_0$, $a=A_*$, $b=B_* r_*^{1/s}$ and $\alpha=1/s$ then we 
can apply Lemma~\ref{lem:iteration}, to conclude:
\begin{equation}
  \label{eq:secondest}
  \rho (z(t;z_0)) \leq 2\rho_0 < C_1 r_{**} \quad \forall |t|\leq T_*=
  \rho_0 A_*^{-1} \,  
exp \Big \{ B_* \left( \frac{r_*}{2\rho_0}\right)^{1/s} \Big \} \, .
\end{equation}
This implies that $\mathcal{H}(z(t;z_0))$ is well defined in this interval of 
time, it is not constant and $|\mathcal{H}(z(t;z_0))| \leq < 3 C_1 r_{**}$. 
Recalling that $\mathcal{H}(z)$ is tangent to the identity close to zero, we have
$|z|\leq C_3|\mathcal{H}(z)|$ for some positive $C_3$.

Then setting $A_{**}=2A_*r_{**}^{-1}$, $B_{**}=B_* 
\left( r_*/(2r_{**})\right)^{1/s}$ and $C_{**}=C_1C_3$, we get:
\begin{equation*}
|z(t;z_0)|\leq C_{**}r_{**} \, ,
\end{equation*}
for all $|t|\leq A_{**}^{-1}exp\Big \{
B_{**}\left(\frac{r_{**}}{|z_0|}\right)^{1/s} \Big \}$.

\section{Arithmetical conditions}
\label{sec:conclusions}

In this paper we proved that any analytic germs of vector field of 
$(\C^n,0)$ with diagonal, non--resonant linear part has an 
{\em effective stability} domain, i.e. stable up to finite but ``long times'',
 close to the stationary point, provided the linear part verifies an 
arithmetical Bruno--like 
condition depending on a parameter $s>0$, which in the case of $2$
dimensional vector fields can be put in the form:
\begin{equation*}
\limsup_{n\rightarrow +\infty} \left( \sum_{j=0}^{k(n)}\frac{\log
  q_{j+1}}{q_{j}} - s\log n \right) <+\infty \, . 
\end{equation*}

\begin{remark}[Invariance of $\bigcup_{s>0}\mathcal{B}_s$, $n=1$ under the action of $PSL(2,\Z)$]
\label{rem:invariance}
The continued fraction development~\cite{HardyWright,MMY} of an
irrational number $\omega$ gives us the sequences: $(a_k)_{k\geq 0}$
and $(\omega_k)_{k\geq 0}$. Then we introduce $(\beta_{k})_{k\geq -1}$
defined by $\beta_{-1}=1$ and for all integer $k\geq 0$:
$\beta_{k}=\prod_{j=0}^k \omega_k$, which verifies :
$1/2<\beta_kq_{k+1}<1$ and $q_n \beta_{n-1}+q_{n-1}\beta_{n}=1$, where
$q_k$'s are the denominators of the 
continued fraction development of $\omega$. We claim that
condition Bruno--$s$~\eqref{eq:brunos1dim} is equivalent to the
following one: 
\begin{equation}
  \label{eq:brunosbeta}
  \limsup_{k \rightarrow +\infty}\left( \sum_{j=0}^k \beta_{j-1} \log
  \omega_j^{-1} + s \log \beta_{k-1} \right) < +\infty \, . 
\end{equation}
This can be proved by using the relations between $\beta_l$ and $q_l$,
to obtain the bound, for all integer $k>0$:
\begin{equation*}
 \Big\lvert \sum_{l=0}^k\left( \beta_{l-1}\log
\omega_l +\frac{\log
  q_{l+1}}{q_l}\right)\Big\rvert \leq
\sum_{l=0}^k\Big\lvert\beta_{l-1}\log{\beta_lq_{l+1}}\Big\rvert+\Big\lvert\beta_{l-1}\log  
\beta_{l-1}\Big\rvert+\Big\lvert\frac{q_{l-1}}{q_l}\beta_l \log
q_{l+1} \Big\rvert\notag \leq 18\, , 
\end{equation*}
where we used the convergence of series $\sum q_l^{-1}$ and $\sum
q_l^{-1}\log q_l$ (see~\cite{MMY} page 272).

To prove the invariance of $\bigcup_s \mathcal{B}_s$ under the action
of $PSL(2,\Z)$, is enough to consider its generators: $T\omega =
\omega +1$ and $S\omega = 1/\omega$. For any irrational $\omega$, $T$
acts trivially being $\beta_k(T\omega)=\beta_k(\omega)$ for all $k$,
whereas for $S$ we have
$\beta_k(\omega)=\omega \beta_{k-1}(S\omega)$ for all $k\geq 1$. Let $\omega$ be an
irrational and let $\omega^{\prime}=\omega^{-1}$, let us also denote
with a $\prime$ quantities given by the continued fraction algorithm
applied to $\omega^{\prime}$, then 
using~\eqref{eq:brunosbeta} one can prove:
\begin{equation*}
 \omega_0 \left( \sum_{j=0}^k \beta_{j-1}^{\prime} \log
  {\omega^{\prime}_j}^{-1} + s\omega_0^{-1} \log \beta_{k-1}^{\prime}
  \right)=C(\omega,s) +
\sum_{j=0}^{k+1} \beta_{j-1} \log
  \omega_j^{-1} + s\log \beta_{k}\, ,
\end{equation*}
where $C(\omega,s)=\omega_0\left(\log \omega_1^{-1}-s\log
  \omega_0\right)+\sum_{l=0}^1\beta_{l-1}\log
  \omega_l^{-1}$, from which the claim follows.
\end{remark}

Let us consider a slightly stronger version of the Bruno--$s$ condition: 
$\omega \in (0,1)\setminus \Q$ belongs to $\Tilde{\mathcal{B}}_s$ if:
\begin{equation}
    \label{eq:newbruno1}
  \lim_{n\rightarrow +\infty}\left( \sum_{l=0}^{k} \frac{\log
  q_{l+1}}{q_l}-s \log q_k \right) < +\infty \, , 
\end{equation}
where $(q_n)_n$ are the convergents to $\omega$, and let us introduce
a second arithmetical condition denoted by $\mathcal{B}_s^{\prime}$ to
be the set of irrational numbers whose convergents verify: 
\begin{equation}
  \label{eq:brunoprimes}
 \lim_{k\rightarrow +\infty} \frac{\log q_{k+1}}{q_k\log q_k}= s \, .
\end{equation}

We state without proof the following proposition, and we refer
to~\cite{Carletti2002}, to all details:
 
\begin{proposition}
  \label{prop:differentbruno1s}
Let $s>0$ and let $\omega \in (0,1)\cap \Tilde{\mathcal{B}}_s$. Then
if $\omega$ is not a Bruno number then $\omega \in
\mathcal{B}_s^{\prime}$, otherwise $\omega \in
\mathcal{B}^{\prime}_0$. 
\end{proposition}

Therefore if $\omega \in \Tilde{\mathcal{B}}_s\setminus \mathcal{B}$ then
the denominators of the convergent to $\omega$ can grow like a factorial,
 more precisely, $q_{k+1}=\mathcal{O}\Big( (q_k!)^s \Big)$, is allowed.

\end{document}